\documentclass[12pt]{amsart}

\usepackage{amssymb}
\usepackage{latexsym}
\usepackage{verbatim}
\usepackage{fullpage}

\input {cyracc.def}
 \font\tencyr=wncyr10 
\font\tencyi=wncyi10 
\font\tencysc=wncysc10 
\def\rus{\tencyr\cyracc}
\def\rusi{\tencyi\cyracc}
\def\rusc{\tencysc\cyracc}

\newcommand{\re}[1]{\textrm  (\ref{#1})}

\renewenvironment{proof}
{\noindent {\sl Proof.}\quad }{\hfill
$\square$ \vskip1.1ex\noindent }

\newenvironment{proof*}
{\noindent {\sl Proof.}\quad }{\hfill
$\square$}

\renewcommand{\theequation}{\thesection .\arabic{equation}}
\renewcommand{\thesubsubsection}{\theequation .\arabic{subsubsection}}

\catcode`\@=\active
\catcode`\@=11
\def\@eqnnum{\hbox to
.01pt{}\rlap{\hskip-\displaywidth(\mathbf{\theequation})}}
\catcode`\@=12

\newenvironment{s}[1]
{ \vskip1.2ex \refstepcounter{equation}
\noindent {\bf \theequation\quad #1.} \begin{sl}}{\end{sl}
\vskip1.1ex\noindent }

\newenvironment{subs}[1]
{\vskip1.2ex \refstepcounter{equation}
\noindent {\bf (\theequation)\quad #1.} }{\quad}

\newcommand {\ah}{{\frak a}}
\newcommand {\be}{{\frak b}}

\newcommand {\g}{{\frak g}}
\newcommand {\h}{{\frak h}}

\newcommand {\p}{{\frak p}}

\newcommand {\te}{{\frak t}}

\newcommand {\esi}{\varepsilon}
\newcommand {\ap}{\alpha}

\newcommand {\lb}{\lambda}

\newcommand {\vp}{\varphi}

\newcommand {\ch}{{\mathcal H}}

\newcommand {\ad}{{\mathrm{ad\,}}}

\newcommand {\hot}{{\mathrm{ht\,}}}

\newcommand {\rk}{{\mathrm{rk\,}}}

\newcommand {\tri}{{\frak sl}_2}
\newcommand {\GR}[2]{{\textrm{{\bf #1}}}_{#2}}

\newcommand {\ov}{\overline}

\newcommand {\Ab}{{\frak Ab}}
\newcommand {\Abo}{\overset{o}{{\frak Ab}}}

\newcommand {\vno}[1]{\vskip#1 ex\noindent}

\newcommand {\qus}{\hfill $\square$ \vno{1.1}}

\newcommand {\beq}{\begin{equation}}
\newcommand {\eeq}{\end{equation}}

\newfam\Bbbfam\newfam\eufam%
\font\Bbbfont=msbm10 scaled 1200%
\font\olala=msam10 scaled 1200%
\font\frak=eufm10 scaled 1400%
\font\Bbbsmallfont=msbm8%
\textfont\Bbbfam=\Bbbfont\scriptfont\Bbbfam=\Bbbsmallfont
\font\euzw=eufm10 scaled 1200%
\font\euac=eufm7 scaled 1200%
\font\euacc=eufm7 scaled 1000%
\textfont\eufam=\euzw\scriptfont\eufam=\euac%
\scriptscriptfont\eufam=\euacc%
\def\frak{\fam\eufam}%
\def\Bbb{\fam\Bbbfam}%

\def\varnothing{\hbox {\Bbbfont\char'077}}
\def\square{\hbox {\olala\char"03}}

\def\cyeq{\hbox {\olala\char'064}}

\begin{document}
\setlength{\parskip}{2pt plus 4pt minus 0pt}
\hfill {\scriptsize November 3, 2002}
\vskip1ex
\vskip1ex

\title[]{Abelian ideals of a Borel subalgebra and long positive roots
}
\author[]{\sc Dmitri I. Panyushev}
\thanks{This research was supported in part by the Alexander von
Humboldt-Stiftung and RFBI Grant no. 01--01--00756}
\maketitle
\begin{center}
{\footnotesize
{\it Independent University of Moscow,
Bol'shoi Vlasevskii per. 11 \\
121002 Moscow, \quad Russia \\ e-mail}: {\tt panyush@mccme.ru }\\
}
\end{center}

Let $\be$ be a Borel subalgebra of a simple Lie algebra $\g$.
Let ${\Ab}$ denote the set of all Abelian ideals of $\be$.
It is easily seen that any $\ah\in \Ab$ is actually contained in the nilpotent
radical of $\be$. Therefore $\ah$ is determined by the the corresponding
set of roots. More precisely, let $\te$ be a Cartan subalgebra of $\g$
lying in $\be$ and let $\Delta$ be the root system of the pair $(\g,\te)$.
Choose $\Delta^+$, the system of positive roots,
so that the roots of $\be$ are positive.
Then $\ah=\oplus_{\gamma\in I}\g_\gamma$, where $I$ is a suitable subset
of $\Delta^+$
and $\g_\gamma$ is the root space for $\gamma\in\Delta^+$.
It follows that there are finitely
many Abelian ideals and that any
question concerning Abelian ideals can be stated in terms of
combinatorics of the root system.

An amazing result of D.\,Peterson says that the cardinality of
$\Ab$ is $2^{\rk\g}$. His approach uses a one-to-one correspondence between
the Abelian ideals and the so-called `minuscule' elements of the affine
Weyl group $\widehat W$.
An exposition of Peterson's results is found in \cite{Ko1}. Peterson's work
appeared to be the point of departure for active recent investigations of
Abelian ideals, $\ad$-nilpotent ideals, and related problems of
representation theory and combinatorics \cite{AKOP},%
\cite{CP},\cite{CP2},\cite{CFP},%
\cite{Ko1},\cite{KOP},\cite{OP},\cite{pr}.
We consider $\Ab$ as poset with respect to inclusion, the zero ideal
being the unique minimal element of $\Ab$.
Our goal is to study this poset structure. It is easily seen that $\Ab$
is a ranked poset; the rank function attaches to an ideal its dimension.
It was shown in \cite{pr} that
there is a one-to-one correspondence between the maximal Abelian ideals
and the long simple roots of $\g$. (For each simple Lie algebra, the
maximal Abelian ideals were determined in \cite{Ro}.)
This correspondence possesses a number of
nice properties, but the very existence of it was demonstrated in a
case-by-case fashion.  Here we give a conceptual explanation of that empirical
observation. More generally, we prove that

$\bullet$ there is a natural mapping $\tau: \Abo\to \Delta^+_l$,
where $\Abo$ is the set of all nontrivial Abelian ideals and $\Delta^+_l$
is the set of long positive roots, see Proposition~\ref{floor-1}. We say
that $\tau(I)$ is the rootlet of $I$;

$\bullet$ Each fibre $\Ab_\mu:=\tau^{-1}(\mu)$ is a poset in its own right,
and we prove that $\Ab_\mu$ contains a unique maximal and a unique minimal
element, see Theorem~\ref{main}.

$\bullet$ If $I$ is a maximal Abelian ideal, then $\tau(I)$ is a (long)
simple root. Restricting $\tau$ to $\Ab_{max}$, the set of
maximal Abelian ideals, yields the above correspondence;
\\[.6ex]
The uniqueness of maximal and minimal elements suggests that they can
have a nice description.
For any $\mu\in\Delta^+_l$, we explicitly describe the minimal ideal in
$\Ab_\mu$ and the corresponding minuscule element of
$\widehat W$ (Theorem~\ref{unique}).
Let $I(\mu)_{min}$ denote the minimal element of
$\Ab_\mu$. The collection of these ideals has a transparent characterisation:
Given $I\in \Abo$, we have $I=I(\mu)_{min}$ for some $\mu$ if
and only if all roots  of $I$ are not orthogonal to $\theta$, the highest root
(see Theorem~\ref{H}).  We also determine the generators of
the ideals $I(\mu)_{min}$.
\par
In Section~\ref{more}, the structure of posets $\Ab_\mu$ is considered.
It is shown that $\#(\Ab_\mu)>1$ if and only if $(\mu,\theta)=0$.
A criterion is also given for $\#(\Ab_\mu)>2$. In fact, I can give a general
description of $\Ab_\mu$ and, in particular, of
the maximal element $I(\mu)_{max}\in \Ab_\mu$. This description
is in accordance with (actually, is inspired by)
my computations for all simple Lie algebras, but
I cannot give yet a general case-free proof. This description shows
that any $\Ab_\mu$ is isomorphic to the poset of all ideals sitting inside
of an Abelian nilpotent radical.
More precisely, there are a
regular\footnote{this means that the subalgebra is normalized by $\te$}
simple subalgebra $\g_{(\mu)}\subset\g$ and a maximal parabolic
subalgebra $\p_{(\mu)}\subset \g_{(\mu)}$ with Abelian nilpotent radical
$\p^{nil}_{(\mu)}$  such that $\Ab_\mu$ is isomorphic to the
poset of all Abelian $\be_{(\mu)}$-ideals
in $\p^{nil}_{(\mu)}$, see Section~\ref{more} for details.
As is well-known,
the latter is isomorphic to the weight poset of a fundamental
representation of the Langlands dual Lie algebra $\g_{(\mu)}^\vee$
\cite{bob},\,\cite{stembr}.
Since this fundamental representation is minuscule, the
weight poset of it is isomorphic to the Bruhat poset
$W^{(\mu)}/W^{(\mu)}_\vp$. Here
$W^{(\mu)}$ is the Weyl group of $\g_{(\mu)}$ (or  $\g_{(\mu)}^\vee$)
and $ W^{(\mu)}_\vp$ is the stabilizer of the fundamental weight in question.
Such posets are also called {\it minuscule\/}.
This completely solves the problem of decribing the structure of $\Ab_\mu$.
\par
In Section~\ref{primery}, the general theory developped so far is illustrated
with examples related to all simple Lie algebras. We compute $\#(\Ab_\mu)$
for each $\mu\in\Delta^+$.
For ${\frak sl}_n$, ${\frak sp}_{2n}$, $\GR{G}{2}$, and $\GR{F}{4}$,
an explicit description of the posets
$\Ab_\mu$ is given. In case  of
${\frak sl}_n$, an algorithm is presented for writing out the minuscule element
corresponding to an Abelian ideal.
\\
Our proofs are based on the relationship between the Abelian ideals and the
minuscule elements in the affine Weyl group. We repeatedly use the
procedure of extension of Abelian ideals that  follows from this relationship.

\section{Preliminaries on Abelian ideals}

\noindent
\begin{subs}{Main notation}
\end{subs}
$\Delta$ is the root system of $(\g,\te)$ and
$W$ is the usual Weyl group. For $\ap\in\Delta$, $\g_\ap$ is the
corresponding root space in $\g$.

$\Delta^+$  is the set of positive
roots and $\rho=\frac{1}{2}\sum_{\ap\in\Delta^+}\ap$.

$\Pi=\{\ap_1,\dots,\ap_p\}$ is the set of simple roots in $\Delta^+$.
 \\
 We set $V:=\te_{\Bbb Q}=\oplus_{i=1}^p{\Bbb Q}\ap_i$ and denote by
$(\ ,\ )$ a $W$-invariant inner product on $V$.   As usual,  $\mu^\vee=2\mu/(\mu,\mu)$ is the coroot
for $\mu\in \Delta$.
Letting $\widehat V=V\oplus {\Bbb Q}\delta\oplus {\Bbb Q}\lb$, we extend
the inner product $(\ ,\ )$ on $\widehat V$ so that $(\delta,V)=(\lb,V)=(\delta,\delta)=
(\lb,\lb)=0$ and $(\delta,\lb)=1$.

$\widehat\Delta=\{\Delta+k\delta \mid k\in {\Bbb Z}\}$ is the set of affine real roots
and $\widehat W$ is the  affine Weyl group.
\\
Then $\widehat\Delta^+= \Delta^+ \cup \{ \Delta +k\delta \mid k\ge 1\}$ is
the set of positive
affine roots and $\widehat \Pi=\Pi\cup\{\ap_0\}$ is the corresponding set
of affine simple roots.
Here $\ap_0=\delta-\theta$, where $\theta$ is the highest root
in $\Delta^+$.  The inner product $(\ ,\ )$ on $\widehat V$ is
$\widehat W$-invariant.
\\
For $\ap_i$ ($0\le i\le p$), we let $s_i$ denote the corresponding simple reflection in $\widehat W$.
If the index of $\ap\in\widehat\Pi$ is not specified, then we merely write
$s_\ap$. 
The length function on $\widehat W$ with respect
to  $s_0,s_1,\dots,s_p$ is denoted by $l$.
For any $w\in\widehat W$, we set
\[
\widehat N(w)=\{\ap\in\widehat\Delta^+ \mid w(\ap) \in -\widehat \Delta^+ \} .
\]
If $w\in W$, then $\widehat N(w)\subset \Delta^+$ and we also write
$N(w)=\widehat N(w)$ in this case.
\begin{subs}{Abelian ideals}
\end{subs}
Let $\ah\subset\be$ be an Abelian ideal. It is easily seen that
$\ah\subset [\be,\be]$. Therefore $\ah=\underset{\ap\in I}{\oplus}\g_\ap$
for a subset $I\subset \Delta^+$, which is called the {\it set of roots of\/}
$\ah$.
As our exposition will be mostly combinatorial, an Abelian ideal
will be identified with the respective set of roots.
That is, $I$ is said to be an Abelian ideal, too.  Whenever we want to
explicitly indicate the context,
we say that $\ah$ is a {\it geometric\/}
Abelian ideal, while $I$ is a {\it combinatorial\/} Abelian ideal.
In the combinatorial context, the definition of an Abelian ideal
(subalgebra) can be stated as follows.
\\
$I\subset\Delta^+$ is an Abelian ideal, if the following two conditions
are satisfied:

(a) for any $\mu,\nu\in I$, we have $\mu+\nu\not\in\Delta$;

(b) if $\gamma\in I$, $\nu\in\Delta^+$, and $\gamma+\nu\in\Delta$, then
$\gamma+\nu\in I$.  \\
If $I$ satisfies only (a), then it is called an Abelian {\it subalgebra\/}.
\\[.6ex]
Following D.\,Peterson, an element $w\in \widehat W$ is said to be
{\it minuscule},
if $\widehat N(w)$ is of the form $\{\delta-\gamma \mid \gamma\in I \}$,
where $I$ is a subset of $\Delta^+$. It was shown by Peterson
that such an $I$ is a combinatorial Abelian
ideal and, conversely, each Abelian ideal occurs in this way,
see \cite[Prop.\,2.8]{CP},\,\cite{Ko1} .
Hence one obtains a one-to-one correspondence between the
Abelian ideals of $\be$ and the minuscule elements of $\widehat W$. If $w\in
\widehat W$ is minuscule, then $I_w$ (resp. $\ah_w$) is the corresponding
combinatorial
(resp. geometric) Abelain ideal. That is,
\[
   I_w=\{\gamma\in\Delta^+ \mid \delta-\gamma \in\widehat N(w)\}
\ \textrm{ and }\  \ah_w=\oplus_{\ap\in I_w}\g_\ap \ .
\]
Conversely, given $I\in \Ab$, we write $w\langle I\rangle$ for the
respective minuscule element.
Notice that
\[
  \dim\ah_w =\# (I_w)=l(w) \ .
\]
Accordingly, being in combinatorial (resp. geometric) context,
we speak about cardinality (resp. dimension) of an ideal.
Throughout the paper, $I$ or $I_w$ stands for a combinatorial Abelian ideal.

\section{Generators of Abelian ideals and long positive roots}
\setcounter{equation}{0}

\noindent
Given an Abelian ideal $I$, let us say that $\gamma\in I$
is a {\it generator}
of $I$, if $\gamma-\alpha\not\in I$ for all $\ap\in\Delta^+$.
Clearly, this is equivalent to the fact that $I\setminus\{\gamma\}$
is still an Abelian ideal. Conversely, if $\varkappa$ is a maximal element
of $\Delta^+\setminus I$ (i.e., $(\varkappa+\Delta^+)\cap \Delta\subset I$)
and $(\varkappa + I)\cap \Delta =\varnothing$,
then $I\cup\{\varkappa\}$ is an Abelian ideal.
These two procedures show that the following is true.

\begin{s}{Proposition}  \label{ranked}
Suppose $I\subset J$ are two Abelian ideals. then there is a chain of
Abelian ideals $I=I_0\subset I_1\subset\ldots\subset I_m=J$ such that
$\#(I_{i+1})=\#(I_i)+1$. In other words, $\Ab$ is a ranked poset, with
cardinality (dimension) of an ideal as the rank function.
\end{s}%
In the geometric setting, the set of generators
has the following description. For an ideal $\ah=\oplus_{\gamma\in I}\g_\gamma
\subset \be$, there is a unique
$\te$-stable space $\tilde\ah\subset \ah$ such that
$\ah=[\be,\ah]\oplus\tilde\ah$. Then $\gamma$ is a generator of
$I$ if and only if it is a weight of $\tilde\ah$.
\\
However, we
need a description of generators of $I$ in terms of the respective
minuscule element. 
As usual, we write $\gamma>0$ (resp. $\gamma<0$), if
$\gamma\in\widehat\Delta^+$ (resp. $\gamma\in -\widehat \Delta^+$).
Let $w\in \widehat W$ be minuscule.
Because $\ap_i\not\in \widehat N(w)$ ($i=1,\dots,p$), any reduced
decomposition of $w$ must end up with $s_0$.
Let $w=s_{i_1}{\cdot}\dots{\cdot}s_{i_r}s_0$ be a reduced decomposition.
As is well known, one then has
\[
\widehat N(w)=\{ \ap_0, s_0(\ap_{i_r}),s_0s_{i_r}(\ap_{i_{r-1}}),\dots,
s_0s_{i_r}\cdots s_{i_2}(\ap_{i_1})\}=
\]
\[
=: \{\delta-\theta,\delta-\gamma_r,\dots, \delta-\gamma_1\}
\]
Here $\gamma_i\in\Delta^+$ and $I_w=\{\theta,\gamma_r,\dots,\gamma_1\}$.
By construction, we have
$\delta-\gamma_1=s_0s_{i_r}\cdots s_{i_2}(\ap_{i_1})$ and hence
$w(\delta-\gamma_1)=-\ap_{i_1}$. Thus,

\parbox[t]{424pt}{\it Any reduced decomposition of $w$ induces a total
ordering on the set $\widehat N(w)$. Moreover, $w$ takes the last element in
$\widehat N(w)$ to $-\widehat\Pi$, i.e., $w(\delta-\gamma_1)=-\alpha_{i_1}$.
}
\\[.7ex]
It follows that if we `shorten' $w$, i.e. consider the element
$w'=s_{i_1}w$, then
$\widehat N(w)=\widehat N(w')\cup\{\delta-\gamma_1\}$ and
$w'(\delta-\gamma_1)=\alpha_{i_1}$. In particular, $w'$ is also minuscule.

\begin{s}{Theorem}  \label{gener}
Suppose $\gamma\in I_w$. Then $\gamma$ is a generator of $I_w$ if and only if
$w(\delta-\gamma)\in-\widehat \Pi$.
\end{s}\begin{proof}
``$\Leftarrow$".
Suppose $w(\delta-\gamma)=-\ap_i$. This means that
$w^{-1}(\ap_i)=\gamma-\delta < 0$. Therefore there exists a reduced
decomposition of $w$ starting with $s_i$: $w=s_iw'$, where $l(w')=l(w)-1$.
Hence $\widehat N(w')=
\widehat N(w)\setminus\{\delta-\gamma\}$ and $w'$ is still a
minuscule element. Thus, $I_w\setminus\{\gamma\}$ is an Abelian ideal.

``$\Rightarrow$".
Suppose $w(\delta-\gamma)\not\in -\widehat\Pi$, i.e.
$w(\delta-\gamma)=-\kappa_1-\kappa_2$, where $\kappa_i\in\widehat\Delta^+$.
Then
$w^{-1}(\kappa_1)+w^{-1}(\kappa_2)=-(\delta-\gamma)< 0$.
Assume for definiteness that $w^{-1}(\kappa_2)< 0$.
Since $w^{-1}(-\kappa_2)>0$ and $w(w^{-1}(-\kappa_2))< 0$, we
have $w^{-1}(-\kappa_2)\in \widehat N(w)$, i.e.
$w^{-1}(-\kappa_2)=\delta-\gamma_2$
for some $\gamma_2\in I_w\subset\Delta^+$. It follows that
$w^{-1}(-\kappa_1)=\delta-\gamma-\delta+\gamma_2\in\Delta$.
As $w(\gamma_2-\gamma)=-\kappa_1 <0$ and $w$ is minuscule,
we must have $\gamma_2-\gamma<0$. Thus $\gamma$ is not a generator of $I_w$.
\end{proof}%

{\sf Remark.} By a result of Cellini and Papi \cite[Theorem\,2.6]{CP}, to
any $\ad$-nilpotent ideal of $\be$ (not necessarily Abelian), one may attach a
unique element of $\widehat W$.
Then one can extend Theorem~\ref{gener} to this setting. However, the
proof becomes
more involved, since the procedure of shortening does not work for
the corresponding elements of $\widehat W$. I hope to consider related
problems in a subsequent publication.

\begin{s}{Theorem}    \label{extension}
Let $I_w$ be an Abelian ideal and $\gamma\in \Delta^+\setminus I_w$. Then
$I_w\cup\{\gamma\}$ is an Abelian ideal if and only if
$w(\delta-\gamma)\in\widehat\Pi$.
\end{s}\begin{proof}
``$\Leftarrow$" \
Suppose $w(\delta-\gamma)=\ap_i$. Then $l(s_iw)=l(w)+1$ and
$\widehat N(s_iw)=\widehat N(w)\cup \{\delta-\gamma\}$.
That is, $s_iw$ is again minuscule
and hence $I_w\cup\{\gamma\}$ is an Abelian ideal.

``$\Rightarrow$" \
It is clear that $\gamma$ is a generator for $I_w\cup\{\gamma\}=:I_{\tilde w}$.
By Theorem~1, we then have $\tilde w(\delta-\gamma)\in -\widehat\Pi$.
Assume that it is $-\ap_i$. Then $w=s_i\tilde w$ and $w(\delta-\gamma)
=\ap_i$.
\end{proof}%
Given a non-trivial minuscule $w\in\widehat W$, it was noticed before that
$w(\ap_i)>0$, $i\in\{1,\dots,p\}$, and  $w(\ap_0)<0$. Let us study the last
element.    Let $\Delta^+_l$ denote the subset of {\it long\/}  roots in $\Delta^+$.
In  the simply-laced case, all roots are proclaimed to be long.

\begin{s}{Proposition}  \label{floor-1}
If $w$ is a non-trivial minuscule element, then $w(\ap_0)+\delta\in \Delta^+_l$.
\end{s}\begin{proof}
Since $w(\ap_0)$ is negative, we can write
$w(\ap_0)=-k\delta-\gamma_0$, where $k\in \{0,1,2,\dots\}$ and
$\gamma_0\in\Delta$.   Recall that $\ap_0=\delta-\theta$.

a) Assume $k\ge 2$. Then
$w(2\delta-\theta)=-(k-1)\delta-\gamma_0<0$. This contradicts the fact
that $w$ is minuscule.

b) Assume $k=0$. Then $w(\delta-\theta)=-\gamma_0$ and
$\gamma_0\in\Delta^+$. It is clear that $w\in \widehat W\setminus W$.
Write the expression of $\theta$ through the simple roots:  $\theta=\sum_{i=1}^p n_i\ap_i$ and set
$\gamma_i=w(\ap_i)$. 
Then $\sum_{i=1}^p n_i\gamma_i=\gamma_0+\delta$. Since $\gamma_i$'s are positive and $\gamma_0\in
\Delta$, there exists a unique  $i_0\in\{1,\dots,p\}$ such that
$n_{i_0}=1$, $\gamma_{i_0}\in\delta+\Delta$ and $\gamma_i\in \Delta$
for $i\ne i_0$.
It follows that the elements $-\gamma_0$, $\gamma_j \ (j\ge 1,\ j\ne i_0)$ form a basis for
$\Delta$. Hence there is $w'\in W$ which takes
$-\gamma_0$, $\gamma_j \ (j\ne i_0)$ to $\ap_1,\dots,\ap_p$.
Because $w'(\gamma_{i_0})\in\delta+\Delta$ and the elements
$w'(\gamma_i)$ $(i=0,1,\dots,p)$ form a basis for $\widehat \Delta$,
we see that $w'(\gamma_{i_0})=\ap_0$.
Thus, $w'w$ takes $\widehat\Pi$ to itself and hence $w'w=1$.
This is however impossible, since $w\not\in W$.

Thus, $k=1$ and $\mu:=w(\ap_0)+\delta=w(2\delta-\theta) \in\Delta$.
Since  $\delta$ is isotropic and $\theta$ is long,  $\mu$ is long as well.
Finally, since $w$ is minuscule, $2\delta-\theta\not\in \widehat N(w)$. Hence $\mu$ is positive.
\end{proof}%
Let $\Abo$ denote the set of all non-trivial Abelian ideals.
By Proposition~\ref{floor-1}, one  obtains the mapping
\[
   \tau : \Abo \to \Delta^+_l \ ,
\]
which is given by
\[
     \tau(I_w)= w(\ap_0)+\delta \ .
\]
The long positive root $\tau(I_w)$ is said to be the {\it rootlet} of the Abelian ideal $I_w$.
Note that the  ideal $\{\theta\}$ is the unique minimal element of $\Abo$
and, by Peterson's result, $\#(\Abo)=2^{\rk\g}-1$.

\begin{s}{Theorem}  \label{not-max}

{\sf 1}. The mapping $\tau$ is onto;

{\sf 2}. If the rootlet of $I_w$ is not simple, i.e.,
$w(\ap_0)+\delta\in\Delta^+\setminus\Pi$,  then $I_w$ is not maximal.

{\sf 3}.
If $\Delta$ is simply-laced and $\tau(I_w)$ is not simple, then there are at
least\/ {\sf two}
maximal Abelian ideals containing $I_w$.
\end{s}\begin{proof}
1. We perform a descending induction on the height of the rootlet of an ideal.
The rootlet with maximal height is $\theta$. Here one takes
$w=s_0$. Then $I_{s_0}=\{ \theta \}$ and $\tau(I_{s_0})=\theta$.
The induction step goes as follows. If $\mu=\tau(I_w)$ and $\mu\not\in \Pi$,
then there exists an $\ap\in\Pi$ such that $(\ap,\mu)>0$.    Then $\mu'=s_\ap(\mu)=\mu-
n_\ap\ap\in \Delta^+_l$ and $\hot (\mu')=\hot(\mu)-n_\ap$. Notice that $n_\ap=1$ if and only if
$\ap$ is long. Set $\mu''=\mu-\ap$. It is again a positive root (not necessarily long).

We have $w(\delta-\theta)=-\delta+\mu''+\ap$. Hence
$w^{-1}(\mu'')+w^{-1}(\ap)=2\delta-\theta$. It follows that

$\left\{
\begin{array}{l}
w^{-1}(\mu'')=(k+2)\delta-\mu_1 \\ 
w^{-1}(\ap)=-k\delta-\mu_2
\end{array} \right. $ \qquad
 for some $k\in\Bbb Z$ and
$\mu_1,\mu_2\in\Delta^+$ such that $\mu_1+\mu_2=\theta$.
\\[.7ex]
As $w$ is minuscule, neither of the elements in the RHS is negative
(for instance, if $w^{-1}(\ap)$ were negative, i.e., $k\ge 0$, then
$w(\mu_2)=-k\delta-\ap_2<0$, which contradicts the fact that $w$ is minuscule).
It follows that $k+2>0$ and $-k>0$, hence $k=-1$.
In particular, we have $w(\delta-\mu_2)=\ap\in\Pi$. It then follows from
Theorem~\ref{extension} that
$w'=s_\ap w$ is again a minuscule element and $I_{w'}=I_w\cup\{\mu_2\}$.
The previous formulae
show that  $\tau(I_{w'})=s_\ap(\mu)=\mu'$. Obviously, any positive long root
can be obtained from $\theta$ through a suitable sequence of
simple reflections. Hence the assertion.

2. The previous argument also shows that if $\tau(I_w)\not\in \Pi$, then
$I_w$ is contained in a larger Abelian ideal.

3. As above, $\mu=\tau(I_w)$.
Making use of the induction argument from part 1, we may reduce the problem to the case, where
$\hot(\mu)=2$. Then $\mu=\ap_1+\ap_2$ -- the sum of two {\it simple} roots. Again the argument from part~1
(with $\ap_1$ and $\ap_2$ in place of $\mu''$ and $\ap$)
shows that there are two different Abelian extensions of $I_w$; namely,
$I_{w_1}=I_w\cup\{\mu_1\}$ and $I_{w_2}=I_w\cup\{\mu_2\}$, where
$w^{-1}(\ap_1)=\delta-\mu_1$ and $w^{-1}(\ap_2)=\delta-\mu_2$.
But $I_w\cup\{\mu_1, \mu_2\}$ is not Abelian, since
$\mu_1+\mu_2=\theta$.
\end{proof}%
{\sf Remark.} In the doubly-laced case, it may happen that
the rootlet of an Abelian ideal is not simple, but the ideal lies in a unique
maximal one. For instance, let $\g$ be the simple Lie algebra of type
$\GR{F}{4}$. We use Vinberg--Onishchik's numbering of simple roots
\cite{VO}. If $\mu=2\ap_2+\ap_3$, then $\tau^{-1}(\mu)$
consists of two ideals (of dimension 7 and 8).
In  the notation of Table~1 in Section~\ref{primery},
$\tau^{-1}(\mu)=\{ I''_7,I'_8\}$. The only maximal ideal containing these two is $I_9$.
\\[.6ex]
Denoting by $\Pi_l$ the set of long simple roots in $\Pi$,
we record an important consequence of the theorem.

\begin{s}{Corollary}   \label{max-ab}
If $I_w$ is a maximal Abelian ideal, then $w(\ap_0)+\delta\in\Pi_l$.
\end{s}%
Thus, denoting by $\Ab_{max}$ the set of all maximal Abelian ideals,
we obtain the mapping
\[
\bar\tau: \Ab_{max}   
\to \Pi_l \ ,
\]
which is the restriction of $\tau$ to $\Ab_{max}$.
By Theorem~\ref{not-max}, $\bar\tau$ is onto.
We shall prove below  that $\bar\tau$ is actually one-to-one.
It turns out that the correspondence obtained between the maximal Abelian
ideals and the long simple roots is precisely the one described in
\cite{pr}. So that our present results provide an {\sl a priori} proof for
some empirical observations in that paper.

\section{Basic properties of posets $\Ab_\mu$}
\label{sect-main}
\setcounter{equation}{0}

\noindent
Given $\mu\in \Delta^+_l$, let $\Ab_\mu$ denote the fibre of $\mu$
for $\tau : \Abo \to \Delta^+_l$. The following useful equality
is a consequence of Peterson's result:
\[
  \sum_{\mu\in\Delta^+_l}\#(\Ab_\mu)=2^{\rk\g}-1  \ .
\]
Each $\Ab_\mu$ is a poset in its own right, and it appears that cutting
$\Abo$ into pieces parametrized by $\Delta^+_l$ has a number of good
properties.

\begin{s}{Theorem}  \label{main}
For any $\mu\in\Delta^+_l$, we have
\begin{itemize}
\item[\sf (i)] the poset $\Ab_\mu$ contains a unique maximal and
a unique minimal element;
\item[\sf (ii)] The dimension of the minimal Abelian ideal in
$\Ab_\mu$ is equal to $1+(\rho,\theta^\vee-\mu^\vee)$.
\item[\sf (iii)] If $I,J\in \Ab_\mu$ and $I\subset J$, then any
intermediate ideal  also belong to $\Ab_\mu$. In particular, $\Ab_\mu$ is a ranked poset.
\end{itemize}
\end{s}%
The proof of this result consists of several parts.
The uniqueness of the minimal (resp. maximal) element will be proved in
Proposition~\ref{unique-min} (resp. Proposition~\ref{unique-max}), and
the dimension formula for the minimal ideal is proved in
Theorem~\ref{unique}. The latter is a by-product of
an explicit description of the minimal ideal
in $\Ab_\mu$ obtained in Section~\ref{describe-min}.
Part (iii) is proved in Corollary~\ref{compare}.
\\[.6ex]
To prove the theorem, we look at the procedure of extension of
Abelian ideals in more details.   If $I,J\in\Ab$, $\dim J=\dim I+1$, and $I\subset J$, then we say
that $J$ is an (Abelian) {\it elementary extension\/} of $I$. Given $I=I_w$, it follows from
Theorem~\ref{extension} that an elementary extension of $I_w$ is possible
if  and only if
$w(\delta-\gamma)=\ap_i\in \widehat\Pi$ for some $\gamma\in \Delta^+$.
Then one can replace $w$ with $w'=s_i w$ and
$I_w$ with $I_{w'}=I_w\cup\{\gamma\}$.
The passage $w \to s_i w$ is also said to be an
{\it elementary extension\/} (via the reflection $s_i$).
Let us realize what happens with the rootlet under this procedure.
Recall that $\Delta$
(or, more generally, the root lattice) has a standard partial order;
one writes $\mu\cyeq \nu$, if $\nu-\mu$ is a sum of positive roots.

\begin{s}{Proposition} \label{elem-ext}
Suppose $I_{w'}$ is an elementary extension of $I_w$, as above.
Then $\tau(I_{w'})=s_i(\tau(I_w))\cyeq
\tau(I_w)$. Moreover, if $w'=s_0w$ (i.e., $i=0$),
then $\tau(I_w)=\tau(I_{w'})$.
\end{s}\begin{proof*}
Set $\nu:=w(\ap_0)+\delta$, the rootlet of $I_w$. Then the rootlet of
$I_{w'}$ is $s_iw(\ap_0)+\delta=s_i(\nu-\delta)+\delta=s_i(\nu)$.
We have two equalities:
$\left\{\begin{array}{l}w(\delta-\gamma)=\ap_i \\ w(\delta-\theta)=\nu-\delta
\end{array}\right.$ . \quad
Consider two possibilities for $i$.

(a) $i\ne 0$. Here we have
\[
  (\ap_i,\nu)=(\ap_i,\nu-\delta)=(\delta-\gamma,\delta-\theta)=
(\gamma,\theta)\ge 0 \ ,
\]
as $\delta$ is isotropic. It follows that
$s_i(\nu)=\nu-(\nu,\ap_i^\vee)\ap_i\,\cyeq\, \nu$.

(b) $i=0$. As $\ap_0=\delta-\theta$, we obtain
\[
  0\le (\gamma,\theta)=(\nu-\delta,\delta-\theta)=-(\nu,\theta)\le 0 \ .
\]
Hence $(\gamma,\theta)=(\nu,\theta)=0$ and
$s_0w(\ap_0)+\delta=s_0(\nu)=\nu$.
\end{proof*}%
\begin{s}{Corollary}    \label{compare}
If $I,J\in\Abo$ and $I\subset J$, then $\tau(J) \cyeq \tau(I)$.
In particular, if $I,J\in \Ab_\mu$, then any intermediate ideal also
belong to $\Ab_\mu$.
\end{s}\begin{proof}
Obviously, for any pair $I\subset J$ of Abelian ideals
there is a sequence of elementary extensions that makes $J$ from $I$.
\end{proof}%
The following result will be our main tool in induction arguments.

\begin{s}{Proposition}  \label{compose}
Let $I=I_w$ be an Abelian ideal.
Suppose $I$ has two different elementary extensions
$I_1=I\cup\{\gamma_1\}$ and $I_2=I\cup\{\gamma_2\}$.
Write $s_iw$ for the minuscule element corresponding to $I_i$, $i=1,2$.

{\sf 1}. If $\tilde I:=I_1\cup I_2$ is not Abelian, then
$\tau(I_1)=\ap_2$, $\tau(I_2)=\ap_1$, and
$\tau(I)=\ap_1+\ap_2$. Moreover, $\ap_1,\ap_2 \in \Pi_l$.

{\sf 2.} If $\tilde I$ is Abelian, then $s_1s_2=s_2s_1$ and  $w\langle \tilde I\rangle=s_1s_2w$;

{\sf 3.} If $\tau(I)=\tau(I_1)$, then $\tilde I$ is Abelian as well and $\tau(I_2)=\tau(\tilde I)$.
\end{s}\begin{proof*}
The equalities $s_iw=w\langle I_i\rangle$ and $I_i=I\cup\{\gamma_i\}$ mean together that
\begin{equation}  \label{dva-ideala}
w(\delta-\gamma_i)=\ap_i\in \widehat\Pi,\ \ i=1,2 \ .
\end{equation}
1. Assume that $I_1\cup I_2$ is not Abelian. Since both $I_1$ and $I_2$ are
Abelian, the only possibility for this
is that $\gamma_1+\gamma_2 \in \Delta^+$.
\begin{itemize}
\item[] If $\gamma_1+\gamma_2\ne\theta$,
then there is an $\ap\in\Pi$ such that $\gamma_1+\gamma_2+\ap$  is a
(positive) root. Then $\gamma_1+\ap\in\Delta$ or
$\gamma_2+\ap\in\Delta$ (Exercise!). If, for instance,
the second condition is satisfied, then
$\gamma_2+\ap\in I$ and $\gamma_1\in I_1$, which contradicts the fact that
$I_1$ is Abelian. Hence $\gamma_1+\gamma_2=\theta$.
\end{itemize}
Now, taking the sum of Equations~\ref{dva-ideala} yields
\[
\ap_1+\ap_2=w(2\delta-\gamma_1-\gamma_2)=w(\delta-\theta)+\delta=\tau (I) \ .
\]
Since $\tau(I)\in\Delta^+_l$, we have $\ap_1,\ap_2\in\Pi_l$.
It follows that $\tau(I_1)=s_1(\ap_1+\ap_2)=\ap_2$ and  $\tau(I_2)=s_2(\ap_1+\ap_2)=\ap_1$.
\\[.6ex]
2. The presence of the
elementary extension $I_1 \mapsto I_1\cup\{\gamma_2\}=\tilde I$ shows that
$w\langle\tilde I\rangle=s_2{\cdot}w\langle I_1\rangle=s_2s_1w$ and
$s_2w(\delta-\gamma_1)\in \widehat\Pi$. The latter means that
$s_2(\ap_1)$ is a simple root. It follows that $s_2(\ap_1)=\ap_1$ and
hence $s_2s_1=s_1s_2$.
\\[.6ex]
3. Under the assumption $\tau(I)=\tau(I_1)$, the first case cannot occur.
Hence $\tilde I$ is Abelian.
Since $s_1,s_2$ commute, we have
$s_2s_1w(\ap_0)+\delta=s_2(\nu)=s_2w(\ap_0)+\delta$, i.e.,
$\tau(\tilde I)=\tau(I_2)$.
\end{proof*}%
\begin{s}{Proposition}  \label{unique-min}
For any $\mu\in \Delta^+_l$, the poset $\Ab_\mu$ has a unique minimal element.
\end{s}\begin{proof}
Assume $\tilde I_1,\tilde I_2$ are two different minimal elements of $\Ab_\mu$.
Clearly $I:=\tilde I_1\cap\tilde I_2$ is again an Abelian ideal, but
$\tau(I)$ is strictly less than $\mu$.
\\
The ideal $\tilde I_1$ can be obtained from $I$ via a chain of elementary
extensions, say
\[
 I\to I\cup\{\varkappa_1\} \to \ldots \to
I\cup\{\varkappa_1,\ldots,\varkappa_n\}=
\tilde I_1 \ .
\]
Similarly, let $I\to I\cup \{\eta_1\}$ be the first step in the chain
of extensions leading from $I$ to $\tilde I_2$.
Set $I(k,0)=I\cup\{\varkappa_1,\ldots,\varkappa_k\}$ and
$I(k,1)=I\cup\{\varkappa_1,\ldots,\varkappa_k,\eta_1\}$, $0\le k\le n$.
By construction, $I(0,1)$ and $I(k,0)$ are Abelian ideals.
Consider the sequence of statements depending on $k$: \\[.7ex]
\hbox to \textwidth{\quad (${\mathcal C}_k$) \hfil
     $I(k,0)\ne \tilde I_1$, \ $I(k,1)$ is Abelian, \ and
$\mu=\tau (\tilde I_1)\cyeq \tau (I(k,1))$. \hfil
}
\vskip.6ex

{\sf Claim.}  For any $k\ge 0$,  (${\mathcal
C}_k$) implies (${\mathcal C}_{k+1}$).
\\[.7ex]
Note that  $({\mathcal C}_0)$ is true. (The last inequality follows from
the equality $\tau(\tilde I_1)=\tau(\tilde I_2)$ and
Corollary~\ref{compare}.) Therefore, granting the claim, we conclude
that $({\mathcal C}_n)$ is also true. But this is nonsense, since
$I(n,0)=\tilde I_1$. This contradiction shows that $\Ab_\mu$ cannot have
two minimal elements. Thus, it remains to prove the Claim.

{\it Proof of the Claim.}  By assumption, we have two elementary extensions:
\[
  I(k,0)\to I(k+1,0)  \textrm{ and } I(k,0) \to I(k,1) \ .
\]
If $w:=w\langle I(k,0)\rangle$, then $w\langle I(k,1)\rangle=s'w$ and $w\langle I(k+1,0)\rangle=s''w$
for some simple reflections $s',s''$.

1. Assume that $I(k+1,1)$ is not Abelian.
Applying Proposition~\ref{compose}(1) to the above triplet of ideals,
we obtain $\tau(I(k,0))=\ap'+\ap''$, $\tau(I(k+1,0))=\ap'$, and $\tau(I(k,1))=\ap''$, where
$\ap',\ap''\in\Pi_l$. Since $I(k+1,0)\subset \tilde I_1$, we
have $\tau (\tilde I_1)=\ap'$. On the other hand, our assumptions
give $\tau(\tilde I_1)  \cyeq \tau (I(k,1))=\ap''$.
Whence $\ap' \cyeq \ap''$. This contradiction shows that $I(k+1,1)$
is Abelian.

2. Since $I(k+1,1)$ is Abelian,  Proposition~\ref{compose}(2) says that $s's''=s''s'$ and
$w\langle I(k+1,1)\rangle=s's''w$. It follows that
\[
  \begin{array}{l} \tau(I(k+1,0))=s''(\tau(I(k,0)))=\tau(I(k,0))-n''\ap'' , \\
  \tau(I(k,1))=s'(\tau(I(k,0)))=\tau(I(k,0))-n'\ap' \
  \end{array}
\]
for some $n',n'' \ge 0$. By the hypothesis,
\[
  \tau(\tilde I_1) \cyeq \tau(I(k,1))=\tau(I(k,0))-n'\ap'  \ ,
\]
and, since $I(k+1,0)\subset \tilde I_1$,
\[
   \tau(\tilde I_1) \cyeq \tau(I(k+1,0))=\tau(I(k,0))-n''\ap''  \ .
\]
Hence   $\tau(\tilde I_1) \cyeq   \tau(I(k,0))-n's'-n''s''=\tau (I(k+1,1))$.

3. If $I(k+1,0)=\tilde I_1$, then the inequalities in the previous part of the proof imply that
\[
   \tau(I(k,0))-n''\ap'' \cyeq \tau(I(k,0))-n'\ap' \ .
\]
Hence $n'=n''=0$.  Then $\mu=\tau(\tilde I_1)=\tau(I(k,0))$. Thus,
$I(k,0)$ is smaller than $\tilde I_1$ and has the same rootlet, which contradicts the minimality
of $\tilde I_1$. Hence $I(k+1,0) \ne \tilde I_1$, and the claim is proved.
\qus
This completes the proof of  Proposition~\ref{unique-min}.
\end{proof}%
In what follows, $I(\mu)_{min}$ stands for the minimal element of $\Ab_\mu$.

\begin{s}{Proposition}   \label{unique-max}  For any $\mu\in\Delta^+_l$,
the poset
$\Ab_\mu$ has a unique maximal element.
\end{s}\begin{proof*}
By Proposition~\ref{unique-min},   any ideal $I\subset \Ab_\mu$
can be obtained from $I(\mu)_{min}$ via a chain of elementary extensions.
Moreover,   it follows from Corollary~\ref{compare} that
each ideal in this chain belong to $\Ab_\mu$.
Another consequence is that if $I,J\in \Ab_\mu$, then $I\cap J\in\Ab_\mu$
as well.

Suppose $I_1,I_2\in\Ab_\mu$. Let us prove that $I_1\cup I_2\in \Ab_\mu$.
Consider the set $I_2\setminus I_1$ and pick there a maximal element
with respect to `$\cyeq$', say $\gamma_2$. Arguing by induction, it suffices to prove that $I_1\cup\{\gamma_2\}$
lies in $\Ab_\mu$. Similarly, take a maximal element
$\nu_1\in I_1\setminus I_2$. Applying Proposition~\ref{compose}(3) to
the ideal $I=I_1\cap I_2\in \Ab_\mu$
and the roots $\nu_1,\gamma_2$, we conclude that
$I\cup\{\nu_1,\gamma_2\}$ is in $\Ab_\mu$. If $I':=I\cup\{\nu_1\}\ne I_1$, then
take a maximal element $\nu_2\in I_1\setminus I'$. Then one applies
Proposition~\ref{compose}(3) to $I'$ and $\nu_2,\gamma_2$.
We eventually obtain $I_1\cup\{\gamma_2\}\in \Ab_\mu$.

Since $I_1\cup I_2\in\Ab_\mu$ for any pair $I_1,I_2\in \Ab_\mu$,
we see that $\Ab_\mu$ has a unique maximal element.
\end{proof*}%
 \begin{s}{Corollary}
The map $\bar\tau :\Ab_{max}\to \Pi_l$ is bijective.
\end{s}\begin{proof}
It follows from Corollary~\ref{max-ab} and Proposition~\ref{unique-max}
that the maximal Abelian ideals are precisely the maximal elements of the
posets $\Ab_\ap$, $\ap\in\Pi_l$.
\end{proof}%
In what follows, $I(\mu)_{max}$ stands for the maximal element of
$\Ab_\mu$. We also say that $I(\mu)_{min}$ is the $\mu$-{\it minimal\/}
and $I(\mu)_{max}$ is the
$\mu$-{\it maximal\/} ideal.

\section{ $\mu$-minimal ideals and their properties}
\label{describe-min}
\setcounter{equation}{0}

\noindent
In this section, an explicit description of
$I(\mu)_{min}$ is given for any $\mu\in \Delta^+_l$.
We also characterise the set of all $\mu$-minimal ideals and find the
generators of $I(\mu)_{min}$.

\begin{s}{Theorem}  \label{uniq}
Let $w\in W$ be an element of minimal length
such that $w(\theta)=\mu$. Then

{\sf 1}. $l(w)=(\rho, \theta^\vee-\mu^\vee)$;

{\sf 2}. $N(w^{-1})=\{\gamma\in\Delta^+\mid (\gamma,\mu^\vee)=-1\}$.
\\
In particular, the set $\{u\in W\mid u(\theta) =\mu\}$ contains
a unique element of minimal length.
\end{s}\begin{proof}
1. Recall that $(\rho,\ap^\vee)=1$ for all $\ap\in\Pi$.
A straightforward calculation shows that, for any $\nu\in\Delta$ and
$\ap\in\Pi$,
\[
  (\rho, s_\ap(\nu)^\vee)= (\rho,\nu^\vee)-(\ap, \nu^\vee) \ .
\]
If $\nu$ is long, then $\vert(\ap, \nu^\vee)\vert\le 1 $.
It follows that, for any $w'\in W$ with the property
$w'(\theta)=\mu$, we have $l(w')\ge (\rho, \theta^\vee-\mu^\vee)$.
On the other hand, if $\mu\in \Delta^+_l$ and $\mu\ne\theta$,
then one can always find an $\ap\in\Pi$ such that $(\ap,\mu^\vee)=1$.
This means that starting with $\mu$ and moving up, one can
reach $\theta$ after applying exactly $(\rho, \theta^\vee-\mu^\vee)$
simple reflections.

2. Set $\Delta^+_\mu(i)=\{\gamma\in\Delta^+\mid (\gamma,\mu^\vee)=i\}$.
We are to show that $\Delta^+_\mu(-1)=N(w^{-1})$. Let us compare the
cardinalities of these two sets. By the first part of the proof,
$\# N(w^{-1})=(\rho,\theta^\vee-\mu^\vee)$. On the other hand,
one has the system of two equations
\vskip.7ex

$\left\{\begin{array}{rcl}
(\rho,\mu^\vee)& = & 1+\frac{1}{2}(\#\Delta^+_\mu(1)-\#\Delta^+_\mu(-1)) \\
2(\rho,\theta^\vee)-2 & = & \#\Delta_\theta(1)=\#\Delta_\mu(1)=
\#\Delta^+_\mu(1)+\#\Delta^+_\mu(-1) \ .
\end{array}
\right.$ \\[.7ex]
The first equality stems from the very definition of $\rho$, whereas in
the second equation we use the fact that $\theta$ is dominant and that
$\mu$ and $\theta$ are $W$-conjugate. From the above system we deduce
that $\#\Delta^+_\mu(-1)=(\rho,\theta^\vee-\mu^\vee)=\# N(w^{-1})$.
\\[.6ex]
On the other hand,
if $\gamma\in \Delta^+_\mu(-1)$, then $(w^{-1}(\gamma),\theta)=-1$.
Hence $w^{-1}(\gamma)$ is negative and $N(w^{-1})\supset \Delta^+_\mu(-1)$.
\end{proof}%
Notice that we also proved that if $u\in W$ is any element taking
$\theta$ to $\mu$, then $N(u^{-1})\supset \Delta^+_\mu(-1)$.
In what follows, we write $w_{\mu}$ for the unique element of
minimal length in $W$ that takes $\theta$ to $\mu$.

\begin{s}{Theorem}   \label{unique}
Set $\tilde w_{\mu}=w_{\mu}s_0\in \widehat W$. Then

{\sf 1}. $\tilde w_{\mu}(\ap_0)+\delta=\mu$;

{\sf 2}. $\tilde w_{\mu}$ is minuscule;

{\sf 3.}  the ideal $I_{\tilde w_\mu}$
is contained  in $\{\gamma \in \Delta^+ \mid (\gamma, \theta)> 0\}$;

{\sf 4.}  $I_{\tilde w_\mu}=I(\mu)_{min}$, the minimal element of $\Ab_\mu$,
and $\# (I_{\tilde w_\mu})=(\rho,\theta^\vee-\mu^\vee)+1$.
\end{s}\begin{proof}
1. Obvious.

2. Suppose $(\rho,\theta^\vee-\mu^\vee)=k\ge 1$
and let $w_{\mu}=s_{i_k}\dots s_{i_1}$
be a reduced decomposition.   We argue by induction on $k$.
Set $u:=s_{i_{k-1}}\dots s_{i_1}\in W$ and $\nu:=u(\theta)$.
Then $l(u)=k-1$ and $s_{i_k}(\nu)=\mu$.
Using Theorem~\ref{uniq}, we obtain
\[
  k-1\ge (\rho,\theta^\vee-\nu^\vee)=(\rho,\theta^\vee-\mu^\vee)-
  (\ap_{i_k},\nu^\vee)=k-(\ap_{i_k},\nu^\vee) \ .
\]
Since $\nu$ is long, $(\ap_{i_k},\nu^\vee)=1$.
It follows  that $(\rho,\theta^\vee-\nu^\vee)=k-1$ and
hence $u=w_\nu$.
Set $\tilde w_\nu=w_\nu s_0$.
By the induction assumption, $\tilde w_\nu$ is minuscule.
To prove that $\tilde w_\mu=s_{i_k}\tilde
w_\nu$ is minuscule, one has to verify
that $\tilde w_\nu(\delta-\gamma_{i_k})=\ap_{i_k}$ for some
$\gamma_{i_k}\in \Delta^+$ (see Theorem~\ref{extension}).
In other words, it should be proved that
$\delta- \tilde w_\nu^{-1}(\ap_{i_k})\in \Delta^+$. As we shall see, this is
a direct consequence of previous formulae.
Indeed, $\tilde w_\nu^{-1}(\ap_{i_k})=s_0 w_\nu^{-1}(\ap_{i_k})$ and
\[
(\theta,w_\nu^{-1}(\ap_{i_k}))=(w_\nu(\theta),\ap_{i_k})=(\nu,\ap_{i_k})>0 \ .
\]
The latter shows that $w_\nu^{-1}(\ap_{i_k})\in\Delta^+$ and
$(\ap_0,w_\nu^{-1}(\ap_{i_k}))<0$. Therefore
$s_0 w_\nu^{-1}(\ap_{i_k})=w_\nu^{-1}(\ap_{i_k})-\theta +\delta$. Thus,
$\delta- \tilde w_\nu^{-1}(\ap_{i_k})=\theta-w_\nu^{-1}(\ap_{i_k})\in
\Delta^+$, and we are done.

3. Again, we argue by induction on $l(w_\mu)$.
Using the notation of the previous part of the proof,
it suffices to observe that
$I_{\tilde w_\mu}=I_{\tilde w_\nu}\cup\{\theta-w_\nu^{-1}(\ap_{i_k})\}$
and $(w_\nu^{-1}(\ap_{i_k}),\theta^\vee)=1$.

4.
If $I_{\tilde w_\mu}$ were not minimal in $\Ab_\mu$, then one could
shorten $\tilde w_\mu$,
so that to obtain a minuscule element giving the ideal with
the same rootlet. But this is impossible for length reason, as $w_{\mu}$ has
minimal possible length among the elements taking $\theta$ to $\mu$. The
dimension of this ideal is already computed in Theorem~\ref{unique}.
Finally, $\# (I_{\tilde w_\mu})=l(\tilde w_\mu)=l(w_\mu)+1$.
\end{proof}%
Thus, we have completed the proof of Theorem~\ref{main}.
\\[.6ex]
Set $\ch=\{ \gamma\in \Delta^+ \mid (\theta, \gamma)>0 \}$.
It is the set of the roots for the standard Heisenberg subalgebra of $\g$.
That is, $\h=\underset{\gamma\in\ch}{\oplus}\g_\gamma$
is a Heisenberg subalgebra of
$\g$. Clearly, $\h$ is a non-Abelian ideal of $\be$.

The previous exposition shows that one has a distinguished collection
of Abelian ideals $\{ I(\mu)_{min} \mid \mu\in \Delta^+_l \}$
and the corresponding subset of minuscule elements of $\widehat W$.
These sets admit the following characterizations:

\begin{s}{Theorem}  \label{H}
The following conditions are equivalent for $I_w\in\Abo$:
\begin{itemize}
\item[\sf (i)] $I_w=I(\mu)_{min}$ for some $\mu\in\Delta^+_l$;
\item[\sf (ii)] $I_w\subset \ch$;
\item[\sf (iii)] $w=w's_0$, where $w'\in W$.
\end{itemize}
\end{s}\begin{proof*}
$(i) \Rightarrow (ii)$.  This is proved in Theorem~\ref{unique}.

$(ii) \Rightarrow (iii)$.  Assume that a reduced decomposition of $w'$ contains $s_0$, say
$w'=w_2 s_0 w_1$. Since   $s_0 w_1 s_0$ is also minuscule (see Section~2), we may assume
without loss that $w_2=1$, i.e., a reduced decomposition of $w'$ begins with $s_0$.
Hence, there is  the elementary extension $w_1 s_0 \to s_0 w_1 s_0$. It was already shown that in this case
one adds to the ideal $I_{w_1 s_0} $ a root which is orthogonal to $\theta$, see
Proposition~\ref{elem-ext}(b).

$(iii) \Rightarrow (i)$.    We argue by induction on $l(w')$. Suppose a reduced decomposition of
$w'$ starts with $s_i$, i.e., $w=s_i w'' s_0$ and $w'' s_0$ is also minuscule.
By the induction hypothesis, $I_{w''s_0}=I(\nu)_{min}$,
 where $\nu=w''(\theta)$. Then
$w'' =w_\nu$ and $l(w'')=(\rho,\theta^\vee -\nu^\vee)$. Set $\mu=s_i(\nu)=
w''(\theta)$.
Then $\mu=\tau(I_w)$ and our goal is to prove that $s_i w''=w_\mu$.
Since $w''s_0 \to s_i w''s_0$ is an elementary extension, we have
$w''s_0(\delta-\gamma)=\ap_i\in \Pi$ for some $\gamma\in \Delta^+$.
It follows that $s_0(\delta-\gamma)\ne \delta-\gamma$. This yields
$(\theta^\vee,\gamma)=1$ and $s_0(\delta-\gamma)=\theta-\gamma$.
Hence $w''(\theta-\gamma)=\ap_i$.
Therefore
\[
 (\ap_i,\nu^\vee)=(w''(\theta-\gamma), w''(\theta^\vee))=
(\theta-\gamma,\theta^\vee)=1 \ .
\]
This equality implies that $(\rho, \theta^\vee-\mu^\vee)=
(\rho, \theta^\vee-\nu^\vee)+1=1+l(w'')=l(w')$.
By Theorem~\ref{uniq}, this means that
$w'=s_iw''\in W$ is the shortest element taking $\theta$ to $\mu$,
and we are done.
\end{proof*}%
\begin{s}{Corollary}
There is a natural one-to-one correspondence between the
Abelian $\be$-ideals in the Heisenberg subalgebra and the long positive
roots. \end{s}%
The next result describes the order relation on the set of
$\mu$-minimal ideals.

\begin{s}{Theorem}  For any $\mu,\nu\in\Delta^+_l$, we have
$I(\mu)_{min}\subset I(\nu)_{min}$ $\Longleftrightarrow$ $\nu\cyeq \mu$.
That is, the poset of $\mu$-minimal elements is
anti-isomorphic to the poset $(\Delta^+_l, \cyeq)$.

\end{s}\begin{proof}
``$\Rightarrow$'' This is contained in Corollary~\ref{compare}.
\\[.6ex]
``$\Leftarrow'$'' Let us show that $w_\nu=w'w_\mu$, where $l(w')=
(\rho,\mu^\vee-\nu^\vee)$. Indeed,

(a) The inequality $l(w')\ge l(w_\nu)-l(w_\mu)=(\rho,\mu^\vee-\nu^\vee)$
is clear.

(b) The opposite inequality can be proved by induction. Set
$\mu-\nu=\sum_{\ap\in\Pi}k_\ap\ap$, where $k_\ap\ge 0$. Since $\vert\mu\vert=
\vert\nu\vert$, we obtain $(\nu, \sum k_\ap\ap)<0$. Hence there exist an
$\ap\in\Pi$ such that $k_\ap>0$ and $(\ap,\nu)<0$.
Then $\nu\cyeq s_\ap(\nu)=\nu+(\vert\mu\vert^2/\vert\ap\vert^2)\ap \cyeq \mu$.
(One should use here the fact that, since $\mu$ and $\nu$ are long,
$k_\ap$ is divisible by $\vert\mu\vert^2/\vert\ap\vert^2$.)
\\[.6ex]
Thus, the minuscule element $\tilde w_\nu$ is obtained from $\tilde w_\nu$
via a sequence of elementary extensions and hence
$I(\mu)_{min}\subset I(\nu)_{min}$.
\end{proof}%
Finally, we give a description of the generators for $\mu$-minimal ideals.
If $w=s_0$, then
$I_{s_0}=\{\theta\}$ and everything is clear. So that we may assume that $\mu\ne\theta$,
i.e., $\tilde w_\mu=w_\mu s_0$ and $w_\mu\ne 1$.

\begin{s}{Proposition}   \label{obraz}
For $\mu\ne\theta$, there is a bijection between the generators of $I(\mu)_{min}$ and
the roots $\ap\in\Pi$ such that $\ap+\mu\in\Delta$ (i.e., $(\ap,\mu^\vee)=-1$).
The generator corresponding to such an $\ap$ is $w_\mu^{-1}(\ap+\mu)$.
\end{s}\begin{proof}
By Theorem~\ref{gener}, $\gamma\in\Delta^+$ is a generator  if and only if
$w_\mu s_0(\delta-\gamma)=-\ap\in\widehat \Pi$. By Theorem~\ref{unique}(3),  $(\gamma,\theta)>0$.
Therefore the LHS is equal to $w_\mu(\theta-\gamma)=\mu-w_\mu(\gamma)$ and
$\mu+\ap=w_\mu(\gamma)\in\Delta$. Hence $\ap\in\Pi$ and $\mu+\ap$ is a root.

This argument can be reversed. Given $\ap\in\Pi$ such that $(\ap,\mu^\vee)=-1$, we set
$\gamma=w_\mu^{-1}(\ap+\mu)$. As $(\ap+\mu,\mu^\vee)\ne -1$, it follows from Theorem~\ref{uniq}(2)
that $\gamma>0$. The rest is clear.
\end{proof}%

\section{More on the structure of $\Ab_\mu$}  \label{more}
\setcounter{equation}{0}

We already know that each $\Ab_\mu$ contains a unique maximal and a unique
minimal element. In this section, we first answer the question:
when is the cardinality of $\Ab_\mu$ equal to 1?
An important observation concerning cardinality stems from
Proposition~\ref{elem-ext}. It was proved there that the elementary
extension via the reflection $s_0$ does not affect the rootlet; and in this
case the rootlet of an ideal has to be orthogonal to $\theta$.
What we prove now is that this gives a necessary and sufficient condition
for $\#(\Ab_\mu)>1$.

\begin{s}{Theorem} \label{cardinal}
{\sf (i)} $\#(\Ab_\mu)>1$ if and only if $(\mu,\theta)=0$
(i.e., $\mu\not\in\ch$).
\\
{\sf (ii)} If $(\mu,\theta)=0$, then the non-empty poset $\Ab_\mu \setminus
\{I(\mu)_{min}\}$ has a unique minimal element, say $I'$.
Here $I'=I(\mu)_{min}\cup\{\gamma\}$, where $\gamma=w_\mu^{-1}(\theta)$.
The corresponding minuscule element is $s_0\tilde w_\mu=s_0w_\mu s_0$.
\end{s}\begin{proof}
(i) In view of Theorem~\ref{main}, it is clear that
$\#(\Ab_\mu)>1$ if and only if $I(\mu)_{min}$ has an elementary extension
that does not change the rootlet. So, we stick to considering possible
elementary extensions of $I(\mu)_{min}$. This is based on the explicit
description in Theorem~\ref{unique}.
\\[.6ex]
1. Since $\mu$ is the rootlet, we have
\begin{equation}   \label{rootlet}
   \tilde w_\mu(\delta-\theta)=\mu-\delta \ .
\end{equation}
Suppose there is an elementary extension of $I(\mu)_{min}$, i.e., we have a
$\gamma\in \Delta^+$ such that
\begin{equation}   \label{gamma-ext}
\tilde w_\mu(\delta-\gamma)=\ap\in\widehat \Pi \ .
\end{equation}
There are two possibilities for $\ap$.

(a)  $\ap=\ap_i\in \Pi$. Rewriting Eq.\,\re{gamma-ext} as
$s_0(\delta-\gamma)=w_\mu^{-1}\ap_i$, we see that
$s_0(\delta-\gamma)\in \Delta$. This can only happen if
$(\ap_0,\delta-\gamma)>0$, i.e., $(\theta,\gamma)>0$
(and then $s_0(\delta-\gamma)=\theta-\gamma$\,). Combining
Equations~\re{rootlet} an \re{gamma-ext}, we obtain
$(\mu,\ap_i)>0$ and hence $s_i(\mu)\ne\mu$. Thus, any elementary extension
via a simple reflection from $W$ changes the rootlet of
$I(\mu)_{min}$.

(b) $\ap=\ap_0$. Here we get the
chain following inequalities:
\[
   0\le (\theta,\gamma)=(\delta-\theta,\delta-\gamma)=
   (\mu-\delta,\delta-\theta)=-(\mu,\theta)\le 0 \ .
\]
Thus, we have the conclusion: if $I(\mu)_{min}$ has an extension that does not
change the rootlet, then this extension uses the reflection $s_0$ and
the condition $(\mu,\theta)=0$ should be satisfied.
This proves the ``only if'' part.
\\[.6ex]
2. Suppose $(\theta,\mu)=0$. We wish to find an elementary extension of
$I(\mu)_{min}$ that does not change the rootlet $\mu$.
Recall that $\tilde w_{\mu}=
w_\mu s_0$. Take $\gamma=w_\mu^{-1}(\theta)$.
From the description of $w_\mu^{-1}$ (see Theorem~\ref{uniq}(2)),
it follows that $\theta\not\in N(w_\mu^{-1})$, i.e., $\gamma\in\Delta^+$.
Furthermore,
$(\gamma,\theta)=(w_\mu(\gamma),w_\mu(\theta))=(\theta,\mu)=0$.

Hence $\tilde w_\mu(\delta-\gamma)=\delta-\theta=\ap_0$ and
$s_0(\mu)=\mu$. Thus, $I(\mu)_{min}\cup\{\gamma\}$ is an Abelian ideal
lying in $\Ab_\mu$.

(ii) This is essentially proved in the previous part of proof, since $s_0$ is the only possible reflection
that can be used for constructing an elementary extension of $I(\mu)_{min}$  with the rootlet $\mu$.
\end{proof}%
{\sf Remark.} We have proved that, for $I(\mu)_{min}$, there is at most
one elementary extension which lies inside $\Ab_\mu$, and, if exists,
this extension always exploits the reflection $s_0$. But if
$I\in\Ab_\mu$ is not minimal, then there can exist an elementary extension
via $s_i$ ($i\ne 0$) that does not change the rootlet.
\\[.6ex]
Now, we accomplish the following step in describing cardinality of $\Ab_\mu$.
That is, a criterion will be given for $\#(\Ab_\mu) > 2$.
We already know that the condition $(\mu,\theta)=0$ is necessary.

\begin{s}{Proposition}  \label{ge3}
Suppose $\mu\in\Delta^+_l$ and $(\mu,\theta)=0$. Then
\begin{center}
$\#(\Ab_\mu) >2 \Longleftrightarrow  \exists
\ap_i\in\Pi \textrm{ such that }
(\ap_i,\theta)>0  \ \& \ (\ap_i,\mu)=0\ .
$
\end{center}
If these conditions are satisfied, then next element of
$\Ab_\mu$ is
\begin{center}
$I''=I(\mu)_{min}\cup\{ w_\mu^{-1}(\theta),w_\mu^{-1}(\theta-\ap_i)\}$.
\end{center}
\end{s}\begin{proof}
In view of Theorem~\ref{cardinal}(ii), it is clear that $\#(\Ab_\mu)>2$
if and only if
$I'=I_{s_0w_\mu s_0}$ has an elementary extension with the same rootlet.
So, we stick to
considering possible extensions of   $I_{s_0w_\mu s_0}$.

``$\Leftarrow$'' \ We show that $s_is_0w_\mu s_0$ is again minuscule and
the corresponding rootlet is again $\mu$. The second condition is satisfied, since
$(\ap_i,\mu)=$ and hence $s_i(\mu)=\mu$. The condition that $s_is_0w_\mu s_0$ is minuscule
is equivalent, in view of Theorem~\ref{extension}, to that
$s_0w_\mu s_0(\delta-\gamma)=\ap_i$ for some $\gamma\in\Delta^+$, i.e.,
$\delta-s_0w_\mu^{-1}s_0(\ap_i)\in\Delta^+$. Using the definition of $w_\mu$ and the assumptions,
the last expression is equal to $w_\mu^{-1}(\theta-\ap_i)$. Since $(\mu, \theta-\ap_i)=0$,
we deduce from Theorem~\ref{uniq}(2) that $\theta-\ap_i\not\in N(w_\mu^{-1})$, that is,
$w_\mu^{-1}(\theta-\ap_i)$ is positive.

``$\Rightarrow$''  \ Suppose
there is an elementary extension of $I_{s_0w_\mu s_0}$ that does not affect $\mu$,
i.e., there is  a $\gamma\in \Delta^+$ such that
\begin{equation}  \label{0-mu-0}
      s_0 w_\mu s_0(\delta-\gamma)=\ap_i
\end{equation}
and $s_i(\mu)=\mu$. Clearly, $i\ne 0$, i.e., $\ap_i\in \Pi$. Since $s_i(\mu)=\mu$,
we have $(\ap_i,\mu)=0$. Thus, it remains to prove that $(\ap_i,\theta)>0$.
If not, then $(\ap_i,\theta)=0$ and hence $s_0(\ap_i)=\ap_i$.
Then Eq.~\re{0-mu-0} can be written as $\delta-\gamma=s_0 w_\mu^{-1}(\ap_i)$.
As $(\theta, w_\mu^{-1}(\ap_i))=(\mu,\ap_i)=0$, the right-hand side is equal
to $w_\mu^{-1}(\ap_i)\in \Delta$. This contradiction proves that
$(\ap_i,\theta)>0$.
\end{proof}%
{\sf Remark.} If $\g\ne {\frak sl}_n$, then there is only one simple root that is not
orthogonal to $\theta$. In any case, this condition is easy to verify
in practice.
\\[.6ex]
Actually, I can give a description of $I(\mu)_{max}$ and
$\Ab_\mu$, which is consistent
with both the previous results and my computaions in
Section~\ref{primery}, but I cannot find a general case-free proof yet.
In order to provide a stronger motivation and more evidences in favour of the
following description, 
let us look again at previous results of this section.
We have proved that
\begin{itemize}
\item if $(\mu,\theta)=0$, then $\Ab_\mu=\{I(\mu)_{min}\}$;
\item if $(\mu,\theta)>0$ and there is no simple roots $\ap\in\Pi$ such that
$(\theta,\ap)>0$ and $(\ap,\mu)=0$, then $\Ab_\mu=\{I(\mu)_{min},I'\}$,
where $I'=I(\mu)_{min}\cup\{\gamma\}$ and $\gamma=w_\mu^{-1}(\theta)$;
\item if $(\mu,\theta)>0$ and $\ap\in\Pi$ satisfies the conditions
$(\theta,\ap)>0$ and $(\ap,\mu)=0$, then one can further extend $I'$
as follows: $I''=I'\cup\{\gamma'\}$, where $\gamma'=w_\mu^{-1}(\theta-\ap)$.
\end{itemize}
These first steps of constructing extensions
show that each time one adds to $I(\mu)_{min}$
some roots that are orthogonal to $\nu$.  Moreover,
the following is true.

\begin{s}{Proposition}  \label{chain}
Suppose $\ap_1,\dots,\ap_t$ is a chain of simple roots such that
$(\theta,\ap_1)>0$, $(\ap_i,\ap_{i+1})>0$ ($i=1,\dots,t-1$), and
$(\theta,\mu)=(\ap_1,\mu)=\ldots=(\ap_t,\mu)=0$. Then
$\#(\Ab_\mu)\ge t+1$. More precisely,
\[
\{ I^{(0)}, I^{(1)},\dots, I^{(t)} \}\subset \Ab_\mu  \ ,
\]
where $I^{(0)}=I(\mu)_{min}$ and
$I^{(i+1)}=I^{(i)}\cup\{w_\mu^{-1}(\theta-\ap_1-\ldots -\ap_{i})\}$.
\end{s}\begin{proof}
Argue by induction on $t$. The induction step is the same as the proof
of Proposition~\ref{ge3}.
\end{proof}%
After this preparations,
I can state a general description of $I(\mu)_{max}$ and $\Ab_\mu$.
Let $\tilde\Gamma$ be the extended Dynkin diagram of $\g$. It has the ``usual''
nodes that  correspond to the roots in $\Pi$ and the ``extra"
node corresponding to $-\theta$. Let us delete from
$\tilde\Gamma$ all nodes such that the corresponding roots are not orthogonal
to $\mu$. The remaining graph can be disconnected. Let $\Gamma_\mu$
denote the connected component of it that contains the node corresponding
to $-\theta$. For instance, if $(\mu,\theta)>0$, then $\Gamma_\mu=\varnothing$.
Clearly, $\Gamma_\mu$ is the Dynkin diagram of a regular simple Lie subalgebra
of $\g$. Call this subalgebra $\g_{(\mu)}$.
If $\ap_1,\dots,\ap_k$ are all simple roots of $\g$ that correspond to
the usual nodes of $\Gamma_\mu$, then $\{\theta,-\ap_1,\dots,-\ap_k\}$ can
be taken as
a set of {simple} roots for $\g_{(\mu)}$, and one can consider the
respective set of {positive} roots.
Let $\be_{(\mu)}$ be the Borel subalgebra corresponding to the chosen set of positive roots,
and let $\be^-_{(\mu)}$ be the opposite Borel subalgebra. With this convention,
let $\p_{(\mu)}\supset \be_{(\mu)}$ be the maximal parabolic subalgebra
of $\g_{(\mu)}$
determined by $\theta$ (i.e. $\theta$ is the only simple root of $\g_{(\mu)}$ that is not a root of the
Levi subalgebra of $\p_{(\mu)}$).
Let $M_\mu$ be the set of roots of $\p_{(\mu)}^{nil}$,
the nilpotent radical of $\p_{(\mu)}$.
It is obvious that the nilpotent radical
constructed in this way is Abelian, i.e.,  for any $\gamma\in M_\mu$
the  coefficient of $\theta$ can be only 1.
Thus,
\[
 M_\mu=\{ \theta-\sum_{i=1}^k c_i\ap_i \mid c_i\ge 0\}\cap \Delta \ .
\]
Notice that $\{\ap_1,\dots,\ap_k\}$ is a proper subset of $\Pi$, since $\mu\ne
0$. Therefore $M_\mu\subset\Delta^+$.
Then the promised description of $I(\mu)_{max}$ is
\[
   I(\mu)_{max}=I(\mu)_{min}\cup  w_\mu^{-1}(M_\mu) \ .
\]
Furthermore, to get an arbitrary (combinatorial) Abelian ideal in $\Ab_\mu$,
one should take any subset
$A\subset M_\mu$ such that the corresponding geometric subspace
$\oplus_{\gamma\in A}\g_\gamma
\subset \p_{(\mu)}^{nil}$ be $\be^-_{(\mu)}$-stable.
It is easily seen that the last condition is satisfied if and only if
$\oplus_{\gamma\in \ov{A}}\g_\gamma
\subset \p_{(\mu)}^{nil}$ is $\be_{(\mu)}$-stable, where
$\ov{A}=M_\mu\setminus A$. In other words,
$A\subset M_\mu$ gives rise to an element of $\Ab_\mu$ if and only if
$\ov{A}$ is a combinatorial $\be_{(\mu)}$-ideal. It follows that
$\Ab_\mu$ is anti-isomorphic to the poset of $\be_{(\mu)}$-ideals in
$\p_{(\mu)}^{nil}$. Since the latter is symmetric, the prefix ``anti''
can be removed. The posets of ideals in an Abelian nilpotent radical are
known  as {\it minuscule\/} posets, see e.g.
\cite{bob}, \cite{stembr}. The minuscule posets have a number of interesting properties;
they are rank-symmetric, rank-unimodal, Gaussian, Sperner, etc., see
\cite{rstan}.
\\[.6ex]
What can we prove in general in this situation?
First, since each root in $M_\mu$ is orthogonal to $\mu$,
we have, by Theorem~\ref{uniq}(2), that $w_\mu^{-1}(M_\mu)\subset \Delta^+$.
Second, using the definition of $I(\mu)_{min}$, it is not hard to prove that
any subset $I(\mu)_{min}\cup w_\mu^{-1}(A)$
is an Abelian {\it subalgebra\/} of $\Delta^+$.
But it is not theoretically clear why all these subsets are ideals in
$\Delta^+$ and why these lie in $\Ab_\mu$.
However, a direct verification shows that this construction gives the correct
description in all cases.

\section{Examples}   \label{primery}
\setcounter{equation}{0}

\begin{figure}[htb]
\begin{center}
\setlength{\unitlength}{0.021in}
\begin{picture}(115,78)(-5,0)
\multiput(90,0)(15,0){2}{\line(0,1){75}}
\multiput(0,60)(0,15){2}{\line(1,0){105}}
\multiput(0,60)(15,0){6}{\line(0,1){15}}
\multiput(90,0)(0,15){4}{\line(1,0){15}}
\put(94,65){{\small $s_0$}}
\put(76,65){{\small $s_{n-1}$}}
\put(61,65){{\small $s_{n-2}$}}
\put(94,50){{\small $s_1$}}
\put(94,35){{\small $s_2$}}
\put(92,5){{\small $s_{i-1}$}}
\put(5,65){{\small $s_j$}}
\put(93,20){$\cdots$}
\put(19,65){$\cdots$}
\put(1.5,5){{\footnotesize $(i,j)$}}
\thinlines
\qbezier[20](0,0),(0,30),(0,60)
\qbezier[20](15,0),(15,30),(15,60)
\qbezier[30](0,0),(45,0),(90,0)
\qbezier[30](0,15),(45,15),(90,15)
\end{picture}
\caption{The filling of a hook}          \label{hook}
\end{center}
\end{figure}
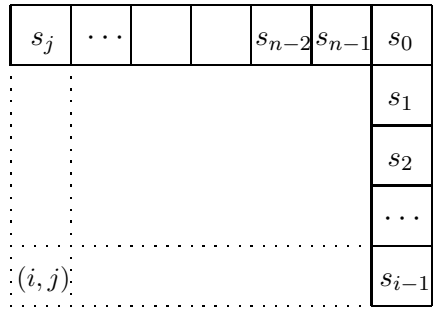

\noindent
Here we present our computations for all simple Lie
algebras.
\begin{subs}{$\g={\frak sl}_n$}
\end{subs}   We assume that $\be$ is the space of
upper-triangular matrices. Then the positive roots are identified
with the pairs $(i,j)$, where $1\le i<j\le n$. Here $\ap_i=(i,i+1)$ and
$\theta=(1,n)$. An Abelian $\be$-ideal is represented by a right-aligned
Ferrers diagram such that the number of rows plus the number of columns is at
most $n$. The unique north-east corner of the diagram corresponds to $\theta$
and the south-west corners give the
generators of the corresponding ideal (see also \cite[3.3]{pr}).
In this case, it is easy to explicitly describe the posets $\Ab_\mu$.
If $\mu=(i,j)$, then
\[
I(i,j)_{max}=\{(p,q) \mid j\le q\ \&\ p\le i\} \textrm{  and  }
I(i,j)_{min}=\{ (1,q) \mid j\le q\}\cup \{ (p,n) \mid 2\le p\le i\} \ .
\]
In other words,
$I(i,j)_{max}$ is the rectangle with the low-left corner at $(i,j)$
and $I(i,j)_{min}$ is the ``north-east'' hook contained in this rectangle, see also Figure~\ref{hook}.
Here $\# I(i,j)_{max}=i(n+1-j)$ and $\# I(i,j)_{min}=n+i-j$.  It follows that
$\# I(i,j)_{max}=\# I(i,j)_{min}$ if and only if  $i=1$ or $j=n$, i.e., precisely
for the roots that are not orthogonal to $\theta$.
It is not hard to compute that
\[
\# \Ab_{(i,j)}= \genfrac{(}{)}{0pt}{}{n+i-j-1}{i-1} \ .
\]
This shows again  that $\#\Ab_{(i,j)}=1$
if and only if $i=1$ or $j=n$. This equality is also in accordance with
Proposition~\ref{ge3}.
It is curious to observe that the assignment $(i,j)\mapsto\#\Ab_{(i,j)}$
gives exactly the Pascal triangle (rotated through the angle $45^o$).
\\[.6ex]
There is an explicit algorithm for writing out the minuscule element for
any $I\in\Abo$. Namely, the minuscule element corresponding to
$I(i,j)_{min}$ is equal to ${(s_{i-1}\ldots s_2s_1)(s_j\ldots
s_{n-2}s_{n-1})s_0}$.
This can be interpreted as a filling of the respective hook,
see Figure~\ref{hook}.

\noindent
Note that the products in parentheses, which correspond to the leg and the
arm of the hook,
commute, so that their order is irrelevant. For an arbitrary Abelian ideal,
one should decompose the corresponding Ferrers diagram as the union of
`north-east' hooks, and then fill in each hook according to the above rule.
The resulting minuscule element is the product of the
corresponding hook elements; the first factor corresponds to the smallest hook, etc.
The best way for understanding all this is to look at the concrete
example.

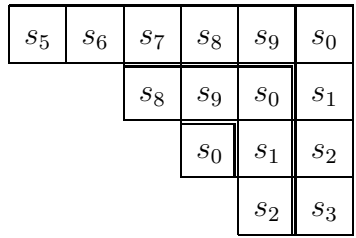
\begin{figure}[htb]
\begin{center}
\setlength{\unitlength}{0.02in}
\begin{picture}(100,60)(-5,0)
\multiput(60,0)(15,0){3}{\line(0,1){60}}
\multiput(0,45)(0,15){2}{\line(1,0){90}}
\multiput(0,45)(15,0){4}{\line(0,1){15}}
\put(30,30){\line(1,0){60}}\put(45,15){\line(1,0){45}}
\put(60,0){\line(1,0){30}}\put(30,30){\line(0,1){15}}
\put(45,15){\line(0,1){30}}
\put(79,50){$s_0$}\put(64,50){$s_9$}\put(49,50){$s_8$}
\put(34,50){$s_7$}\put(19,50){$s_6$}\put(4,50){$s_5$}
\put(79,35){$s_1$}
\put(64,35){$s_0$}
\put(49,35){$s_9$}
\put(34,35){$s_8$}
\put(79,20){$s_2$}
\put(64,20){$s_1$}
\put(49,20){$s_0$}
\put(79,5){$s_3$}
\put(64,5){$s_2$}
\put(30,44){\line(1,0){44}}
\put(74,0){\line(0,1){44}}
\put(45,29){\line(1,0){14}}
\put(59,15){\line(0,1){14}}
\end{picture}
\caption{The decomposition and filling of the Ferrers diagram for an
Abelian ideal in ${\frak sl}_{10}$}  \label{sample}
\end{center}
\end{figure}

Consider the Abelian ideal $I$ in ${\frak sl}_{10}$ with generators
$(1,5)$,\,$(2,7)$,\,$(3,8)$,\,$(4,9)$. Here the Ferrers diagram is
decomposed as the union of three hooks and the corresponding
filling is depicted in Figure~\ref{sample}.
The minuscule element $w\langle I\rangle$ is
$s_0(s_2s_1)(s_8s_9)s_0(s_3s_2s_1)(s_5s_6s_7s_8s_9)s_0$.

\noindent

\begin{subs}{$\g={\frak so}_{2n+1}$ or ${\frak so}_{2n}$}\end{subs}
In the standard notation, the set of long positive roots is
\[
 \Delta^+_l=
\{ \esi_i\pm\esi_j \mid 1\le i< j\le n\}\ .
\]
Here $\theta=\esi_1+\esi_2$ and $\ch\cap\Delta^+_l=\{\esi_i\pm\esi_j\mid
i=1,2 \ \&\ j\ge 3\}
\cup \{\theta\}$.
By Theorem~\ref{cardinal}, $\#(\Ab_\mu)=1$ for any $\mu\in \ch\cap\Delta^+_l$.
By Proposition~\ref{ge3}, we obtain $\# (\Ab_{\esi_1-\esi_2})=2$
and $\# (\Ab_{\esi_3\pm\esi_j})=2$ ($j\ge 4$).
Straightforward computations for the other roots show that
$\# (\Ab_{\esi_i\pm\esi_j})=2^{i-2}$, if $i\ge 3$.
Let us demonstrate how all this is related to the description of
$\Ab_\mu$ in Section~\ref{more}.
\\
Take, for instance,
$\mu=\ap_{n-2}=\esi_{n-2}{-}\esi_{n-1}$ for ${\frak so}_{2n}$.
Then $\g_{(\mu)}=\left\{
\begin{array}{rl}
0, & \textrm{ if } n=4 \\
\tri, & \textrm{ if } n=5 \\
{\frak so}_{2n-6}, & \textrm{ if } n\ge 6 \ .
\end{array}
\right.$
For $n\ge 6$, the Abelian nilpotent
radical in $\g_{(\mu)}$ corresponding to $\theta$ has dimension
$(n-3)(n-4)/2$.
This number is just the difference
$\dim I(\ap_{n-2})_{max}-\dim I(\ap_{n-2})_{min}$. Hence
$\dim I(\ap_{n-2})_{max}=\frac{(n-3)(n-4)}{2} + 2n-3= (n^2-3n+6)/2$, cf.
\cite[Figure\,3]{pr}.
In this case, $\g_{(\mu)}\simeq \g_{(\mu)}^\vee$ and $\#(\Ab_\mu)$ is the
dimension of the half-spinor
representation of ${\frak so}_{2n-6}$, i.e., $2^{n-4}$.

\begin{subs}{$\g={\frak sp}_{2n}$}\end{subs}
In this case, there is only a few long roots:
\[
 \Delta^+_l=\{ 2\esi_i \mid 1\le i\le n\}\ ,
\]
and $\theta=2\esi_1$.
We have $I(2\esi_i)_{min}=\{\esi_1+\esi_i,\dots,
\esi_1+\esi_2, 2\esi_1\}$ and
$I(2\esi_i)_{max}=\{ \esi_k+\esi_j \mid k\le
j\le i\}$. The sole generator of $I(2\esi_i)_{min}$ (resp. $I(2\esi_i)_{max}$)
is $\esi_1+\esi_i$ (resp. $2\esi_i$).
The minuscule element $w\langle I(2\esi_i)_{min}\rangle$ is
$s_{i-1}\ldots s_2s_1s_0$.
Using the matrix presentation of ${\frak sp}_{2n}$ (see e.g. \cite[3.3]{pr}),
it is easily seen that there is a one-to-one correspondence between
the ideals in $\Ab_{2\esi_i}$ and the Abelian ideals of
${\frak sp}_{2i-2}$. Therefore  $\#(\Ab_{2\esi_i})=2^{i-1}$. It is also
possible to give an algorithm for writing out the minuscule element
corresponding to an Abelian ideal in terms of filling of a shifted
Ferrers diagram.

\begin{subs}{$\g=\GR{F}{4}$}\end{subs}
Here we have 12 long positive roots and 15 non-trivial
Abelian ideals. The set $\ch\cap\Delta^+_l$ consists of 9 roots.
Hence the fibre $\Ab_\mu$ contains a unique ideal for these 9 roots
and consists of two ideals for the other 3 roots. The computations of rootlets
and minuscule elements are presented in Table~1.
We follow the numbering of simple roots  from \cite[Tables]{VO},
and the root $\sum_{i=1}^4 c_i\ap_i$ is denoted by $(c_1\,c_2\,c_3\,c_4)$.
For instance, $\theta=(2432)$.
The notation $I_n$ means that the ideal has cardinality $n$.
To distinguish different ideals with the same  cardinality, we use `prime'.
The third, fourth, and fifth columns represent
the ideal, the corresponding minuscule element, and the rootlet, respectively.

\begin{table}[htb]
\begin{center}
\caption{ 
The Abelian $\be$-ideals in $\GR{F}{4}$}
\begin{tabular}{cc|c|c|c}
No.
& $\# I$ & $I$ & $w\langle I\rangle$ & $\tau(I)$ \\ \hline\hline
1\  & 1 & $\{\theta\}$ & $s_0$              & $\theta$ \\
2 & 2 & $\{\theta, 2431\}$ & $s_4s_0$           & $2431$ \\
3 & 3 & $\{\theta, 2431, 2421\}$ & $s_3s_4s_0$        & $2421$ \\
4 & 4 & $\{\theta, 2431, 2421, 2321\}$ & $s_2s_3s_4s_0$     & $2221$ \\ \hline
5 & 5 & $I'_5=I_4 \cup\{2221\}$ & $s_3s_2s_3s_4s_0$ & $2211$ \\
6 & 5 & $I''_5=I_4\cup\{1321\}$ & $s_1s_2s_3s_4s_0$ & $0221$ \\ \hline
7 & 6 & $I'_6=I'_5\cup\{2211\}$ &
    $w'_6=s_4s_3s_2s_3s_4s_0$ & $2210$ \\

8 & 6 & $I''_6=I'_5\cup\{1321\}=I''_5\cup\{2221\}$ &
    $w''_6=s_1s_3s_2s_3s_4s_0$ & $0211$ \\ \hline
9 & 7 & $I'_7 =  I'_6\cup\{2210\}$
    & $w'_7=s_0w'_6$ & $2210$ \\
10 & 7 & $I''_7= I'_6\cup\{1321\}=I''_6\cup\{2211\}$
    & $w''_7=s_1w'_6=s_4w''_6$ & $0210$ \\
11 & 7 & $I'''_7=I''_6\cup\{1221\}$
    & $w'''_7=s_2w''_6$ & $0011$ \\ \hline
12 & 8 & $I'_8=I'_7\cup\{1321\}=I''_7\cup\{2210\}$
     & $w'_8=s_1w'_7=s_0w''_7$ & $0210$ \\
13 & 8 & $I''_8=I''_7\cup\{1221\}=I'''_7\cup\{2211\}$
     & $w''_8=s_2w''_7=s_4w'''_7$ & $0010$ \\
14 & 8 & $I'''_8=I'''_7\cup\{0221\}$
     & $w'''_8=s_3w'''_7$ & $0001$ \\ \hline
15 & 9 & $I_9=I'_8\cup\{1221\}=I''_8\cup\{2210\}$
     & $w_9=s_2w'_8=s_0w''_8$ & $0010$
\end{tabular}
\end{center}
\end{table}

The maximal Abelian ideals are $I'''_8$ and $I_9$.

\begin{subs}{$\g=\GR{G}{2}$}\end{subs}
Here $\#(\Abo)=\#(\Delta^+_l)=3$, so that everything
is easy. Let $\ap$ (resp. $\beta$) be the short (resp. long) simple root. Then
\begin{center}
$\begin{array}{lll}
I_1=\{3\ap+2\beta\}, & w\langle I_1\rangle=s_0, & \tau(I_1)=3\ap+2\beta; \\
I_2=\{3\ap+2\beta,3\ap+\beta\}, & w\langle I_2\rangle=s_\beta s_0, &
\tau(I_2)=3\ap+\beta; \\
I_3=\{3\ap+2\beta,3\ap+\beta,\beta\}, &
w\langle I_3\rangle=s_\ap s_\beta s_0, & \tau(I_3)=\beta;
\end{array}$
\end{center}
\begin{subs}{$\g=\GR{E}{n}$, $n=6,7,8$}\end{subs}
Set $\Delta^+_{(i)}=\{\mu\in\Delta^+ \mid \#(\Ab_\mu)=i\}$ and
$m_i=\#\Delta^+_{(i)}$. Note that $\Delta^+_{(1)}=\ch$. The output of our
calculations of numbers $m_i$ is given in
Table~\ref{m_i}, where we include only the
columns containing nonzero entries. The rightmost column is the control one.

\begin{table}[htb]
\begin{center}
\caption{ 
}  \label{m_i}
\begin{tabular}{c|cccccccc|c}
 & $m_1$ & $m_2$ & $m_3$ & $m_4$ & $m_5$ & $m_6$ & $m_8$ & $m_{12}$ & $\sum
im_i$ \\ \hline
$\GR{E}{6}$ & 21 &  9 &  4 & -- & -- & 2 & -- & -- & $2^6-1$ \\
$\GR{E}{7}$ & 33 & 15 &  8 & 4  & -- & 2 & -- &  1 & $2^7-1$ \\
$\GR{E}{8}$ & 57 & 27 & 16 & 10 &  6 & 3 &  1 & -- & $2^8-1$ \\ \hline
\end{tabular}
\end{center}
\end{table}
%
\noindent
An explicit description of the subsets $\Delta^+_{(i)}$'s
is also obtained. Again, we follow the
numbering of simple roots from \cite{VO}
and denote the root
$\sum_{i=1}^nc_i\ap_i$ by $(c_1\,c_2\ldots c_n)$. For instance, the highest
root of $\GR{E}{6}$ (resp. $\GR{E}{7}$) is $(1\,2\,3\,2\,1\,2)$
(resp. $(1\,2\,3\,4\,3\,2\,2)$).

\begin{table}[htb]
\begin{center}
\caption{ 
}  \label{delta(i)}
\begin{tabular}{c|cccccccc}
 & $\Delta^+_{(1)}$ & $\Delta^+_{(2)}$ & $\Delta^+_{(3)}$ & $\Delta^+_{(4)}$ &
$\Delta^+_{(5)}$ & $\Delta^+_{(6)}$ & $\Delta^+_{(8)}$ & $\Delta^+_{(12)}$  \\
\hline \hline
$\GR{E}{6}$ & $c_6{>}0$ & $\begin{array}{c} c_6{=}0 \\ c_3{>}0
                         \end{array}$ &
$\begin{array}{c}
\{\ap_1{+}\ap_2,\ap_2, \\ \ap_4{+}\ap_5,\ap_4\}
\end{array}$
 & -- & -- & $\{\ap_1,\ap_5\}$ & -- & --  \\ \hline
$\GR{E}{7}$ & $c_6{>}0$ & $\begin{array}{c} c_6{=}0 \\ c_5{>}0
                         \end{array}$ &
$\begin{array}{c} c_6{=}c_5{=}0 \\ c_4>0
\end{array}$
& $\begin{array}{c} c_6{=}c_5{=}c_4{=}0, \\ c_7{>}0\textrm{ or }c_3{>}0
\end{array}$
& -- & $\{\!\ap_2,\ap_1{+}\ap_2\!\}$ & -- &  $\{\ap_1\}$  \\ \hline
$\GR{E}{8}$ & $c_1{>}0$ & $\begin{array}{c} c_1{=}0 \\ c_2{>}0
                         \end{array}$ &
$\begin{array}{c} c_1{=}c_2{=}0 \\ c_3>0
\end{array}$ &
$\begin{array}{c} c_1{=}c_2{=}c_3{=}0, \\ c_4>0
\end{array}$   &
$\begin{array}{l} c_1{=}c_2{=} \\c_3{=}c_4{=}0, \\ c_5>0
\end{array}$
& $\begin{array}{c}\{\ap_8,\ap_6,\\ \ap_6{+}\ap_7\}\end{array}$
& $\{\ap_7\}$ & --  \\ \hline
\end{tabular}
\end{center}
\end{table}

\end{document}